\input amstex
\documentstyle{amsppt}
\magnification 1000
\pagewidth{5.5 true in}
\pageheight{7.5 true in}
\hcorrection{.5 true in}
\vcorrection{.5 true in}
\overfullrule=0pt
\parskip=3pt
\addto\tenpoint{\normalbaselineskip17pt\normalbaselines}
\nologo

\def\t{\roman{Tor}}

\def\l{\lambda}
\def\i{\roman{im}\,}
\topmatter
\title{Bivariate Hilbert Functions for the Torsion Functor}\endtitle
\author{Emanoil Theodorescu}\endauthor
\address{Department of Mathematics, University of Missouri, Columbia, MO, 
65211}
\endaddress
\email{theodore\@math.missouri.edu}\endemail

\abstract{Let $(R,P)$ be a commutative, local Noetherian ring, $I$, $J$ 
ideals, $M$ and $N$ finitely generated $R$-modules. Suppose 
$J + ann_R M + ann_R N$ is $P$-primary. The main result of this paper is 
Theorem 6, which gives necessary and sufficient conditions for the length of 
$\t_i(M/I^nM,N/J^mN)$, to agree with a polynomial, for $m$, $n \gg 0$. As a 
corollary, it is shown that the length of $\t_i(M/I^nM,N/I^nN))$ always 
agrees with a polynomial in $n$, for $n \gg 0$, provided 
$I + ann_R M + ann_R N$ is $P$-primary.}
\endabstract
\endtopmatter

\document
\head{Introduction}
\endhead
\TagsOnRight
 Throughout this paper, unless otherwise stated, $(R,P)$ is a commutative, 
Noetherian local ring with unit and $I$, $J$ are (proper) ideals. Also, let 
$M$, $N$ be finite $R$-modules, $m$, $n$ be nonnegative integers, and let 
$\lambda$ denote length. We would like to study the two-variable  
the Hilbert function $H(n,m):=\l(\t_i(M/I^nM,N/J^mN))$. On the one hand, 
we have in mind extending results on $H(n,m) $ of the authors 
of \cite{BF}, \cite{KS} 
and \cite{WCB}, while on the other hand we seek two variable analogues 
of recent results concerning the Hilbert function 
$H(n) :=\l(\t_i(M/I^nM,N))$. Previous work on $H(n)$ appears in 
\cite{TM}, \cite{VK}) and \cite{ET}. In fact, in \cite{ET} it is shown that 
$H(n)$ agrees with a polynomial in $n$ for $n$ large, if we simply 
assume that the lengths $\l((\t_i(M/I^nM,N))$ are finite. Here we seek to give 
conditions under which $H(n,m)$ has polynomial growth for $n$ and $m$ 
sufficiently large. In some special cases, we give a 
degree bound on the resulting polynomials in $n$ and $m$. 
Determining the exact degree of these polynomials seems to be a 
more difficult task. In the one variable case, \cite {VK} and \cite{ET} give 
upper bound estimates for the degree in 
general and while \cite{ET}, \cite {DK} and \cite {TM} determine the 
degree in some special cases.
 
 In his Doctoral Thesis, Bruce Fields \cite {BF} investigates two-variable  
functions of the form $\l(\t_i(R/I^n, R/J^m))$, where $i \geq 0$, under the 
assumption that $I+J$ is $P$-primary. 
For $i \geq 2$, he proves that these lengths are eventually given by 
polynomials in two variables. Actually, since $\t_i(R/I^n, R/J^m) = 
\t_{i-1}(I^n, R/J^m) = \t_{i-2}(I^n, J^m)$ (by applying twice the shifting 
formula), his proof essentially shows that 
$\oplus_{m,n=0}^{\infty} \t_j(I^nM,J^mN)$, $j \geq 0$, is a finite, bigraded 
module, over a suitable polynomial ring over $R$, where $M$, $N$ are two 
finite $R$-modules. It is then well-known that, if the lengths of homogeneous 
pieces of a finite bigraded module (over a suitable polynomial ring) are 
finite, then they are eventually given by a polynomial function (also see 
Notations and Conventions).

 For $i=0$ and $i=1$, Fields only proves that polynomial growth holds under 
some rather restrictive conditions: he assumes that $R$ is regular local, 
and that $\oplus_{m, n=0}^{\infty} (I^n \cap J^m$) is a {\it finite} bigraded 
module over some polynomial ring in two sets of variables. This is, in 
general, a very strong condition on two ideals $I$, $J$. The function 
$\l(R/(I^n+J^m))$ has also been studied by Kishor Shah \cite{KS} and 
William C. Brown \cite{WCB}, who give sufficient conditions for it to be 
given by a polynomial, for $m$, $n$ $\gg 0$.

The present paper gives a characterization of those cases for which the length 
of \linebreak 
$\t_i(M/I^nM, N/J^mN)$ has polynomial growth, provided the following 
condition is satisfied: $J + ann_R M + ann_R N$ is $P$-primary (see Theorem 6).
It turns out that polynomial growth doesn't always hold, even in the case 
$i \geq 2$, as Fields' work might have suggested (see the Remark following
Corollary 8). On the other hand, Proposition 3 shows that, provided 
$\t_i(I^nM, N/J^mN)$ has finite length, for all large $m$, $n$, its
length is always given by a polynomial, without any restrictive assumption.

 As a corollary to the proof of Theorem 6, under the assumption that 
$I + ann_R M + ann_R N$ is $P$-primary, we prove that 
$\l(\t_i(M/I^nM, N/I^nN))$ has {\it always} polynomial growth. Corollary 8 
shows that, under the hypothesis that both $I+ann_R M+ann_R N$ and 
$J+ann_R M+ann_R N$ be $P$-primary, the length of $\t_i(M/I^nM, N/J^mN)$ has 
polynomial growth if and only if both $\t_i(M, N)$ and $\t_{i-1}(M, N)$ have 
finite length. Finally, when $M \otimes N$ has finite length, Theorem 9 gives 
the formula $\l(\t_i(M/I^nM, N/J^mN)) =$
$$
\l(\t_i(M, N))+\l(\t_{i-1}(I^nM, N))+
\l(\t_{i-1}(M, J^mN))+\l(\t_{i-2}(I^nM, J^mN)), 
$$
\noindent
which works for all $i \geq 0$, by assuming that all $\t_i$ with $i < 0$ are 
zero. 

 The main result of this paper shows that, at least when $J+ann_R M+ann_R N$ 
is $P$-primary, the nature of $\l(\t_i(M/I^nM, N/J^mN))$ is controlled by 
modules of the form $I^nA \cap J^mB$. Therefore, a study of modules of this 
kind would deepen our understanding of $\l(\t_i(M/I^nM, N/J^mN))$.

\head{Notation and Conventions}
\endhead

 We will be using (free) resolutions of modules over several different rings. 
 There will be resolutions of modules over $R$, graded resolutions of graded 
modules over the polynomial ring in $r$ variables, $S_1 := R[X_1, ..., X_r]$, 
as well as bigraded resolutions of bigraded modules over the polynomial ring 
in two sets of variables, $S_2 := R[X_1, ..., X_r; Y_1, ..., Y_s]$. Unless 
otherwise stated, the $\t$'s are over $R$.

 To further simplify notation, we denote $\Cal{M}=\oplus_{n=0}^{\infty} M$, 
which is an (infinitely generated) graded module over the Rees ring 
$\Cal{R}_I:= \oplus_{n=0}^{\infty} I^n$. If $I$ is generated by $x_1$,..., 
$x_r$, then $\Cal{M}$ is naturally an infinitely generated $S_1$-graded 
module, via the canonical ring homomorphism $S_1 @>>> \Cal{R}_I$, given by 
$X_i\mapsto x_i$, for all $i$. The action of $S_1$ on $\Cal{M}$ is given by 
$X_i v_k  = x_iv_k$, where $v_k$ denotes a homogeneous vector of 
degree $k$. Also, if we denote $\Cal{I}\Cal{M} := 
\oplus_{n=0}^{\infty} I^nM$, then this is a finitely generated 
graded module over $\Cal{R}_I$, and hence over $S_1$, as before. It follows 
that $\Cal{M}/\Cal{I}\Cal{M} = \oplus_{n=0}^{\infty} (M/I^nM)$ is a graded 
module over both $\Cal{R}_I$ and $S_1$.

 Similarly, if we assume $J=(y_1, ..., y_s)$, 
$\oplus_{m,n=0}^{\infty} I^nJ^m M$ is a bigraded module over the bigraded 
Rees ring $\Cal{R}_{I,J} := \oplus_{m,n=0}^{\infty} I^nJ^m$, and hence over 
the polynomial ring $S_2$, via a similar map $S_2 @>>> \Cal{R}_{I,J}$.

 Note that any graded free resolution over $S_1$ or $S_2$ of some graded 
module, is also a free resolution of that module over $R$.

 We will be making use of the fact that, in a (bi)graded resolution 
of some $S_1$ (or $S_2$)-graded module, say $\Cal{I}\Cal{M}$, by considering 
just its homogeneous part of degree $k$, we obtain a free resolution, over 
$R$, of the module $I^kM$, the $k$-th homogeneous component of 
$\Cal{I}\Cal{M}$ . 

We will be making repeated use of the fact that, if 
$\Cal{P}:=\oplus_{m, n=0}^{\infty} P_{m, n}$ is a finite 
bigraded $S_2$-module, whose homogeneous pieces have finite length, 
then $\l(P_{m, n})$ is eventually given by a polynomial. In particular, 
$\l(\t_i(I^nM, J^mN))$ is eventually given by a polynomial. Indeed, we can 
take $\Cal{C}$ a $S_1$-graded free resolution (consisting of finite free 
$S_1$-modules) of $\oplus_{n=0}^{\infty} I^nM$ and, similarly, $\Cal{D}$ a 
$S_1^\prime$-graded free resolution of $\oplus_{m=0}^{\infty} J^mN$, also 
consisting of finite free $S_1^\prime$-modules. 
(Here, $S_1^\prime = R[Y_1, ..., Y_s]$.)
Then the modules in $\Cal{C}\otimes_R \Cal{D}$ have a natural structure of 
$S_1\otimes_R S_1^\prime \cong S_2$-modules. Actually, 
$\Cal{C}\otimes_R \Cal{D}$ is a complex of finite, free, $S_2$-modules, whose 
i-th homology is 
$\t_i^R(\oplus_{n=0}^{\infty}I^nM, \oplus_{m=0}^{\infty} J^mN)$. 
Of course, this is a finitely generated bigraded $S_2$-module. 
Since the homogeneous components of this are just $\t_i^R(I^nM, J^mN)$, it 
follows that, if their lengths are finite, then these lengths are eventually 
given by a polynomial in $m$, $n$.

\head{The main result}
\endhead
 In an attempt to study the length of $\t_i(M/I^nM, N/J^mN)$ in as great 
generality as possible, we first investigate $\t_i(I^nM, N/J^mN)$. It 
turns out that in this case polynomial growth follows from the simplest 
assumption that these $Tor$'s have finite length.
The following few results are essentially given without proof, as their proofs 
parallel those of corresponding one-variable statements (see \cite{ET}).

\proclaim{Proposition 1} Let $R$ be a Noetherian ring \rom{(}not necessarily 
local\rom{)}, and $J \subset R$ an ideal. Let $S_1$ be the polynomial ring 
over $R$ in $r$ variables, and let 
$$
\Cal{C}: \qquad \qquad {\Cal F}_2 @>\psi>> {\Cal F}_1 @>\phi>> {\Cal F}_0
$$
\noindent
be a graded complex of graded $S_1$-modules, graded by total degree. Assume 
that ${\Cal F}_1$, ${\Cal F}_0$ are finitely generated $S_1$-modules. 
Then, there is $l \geq 0$, such that, for all $m\geq l$

$$
     H_1(\Cal{C} \otimes \frac{R}{J^m})=\frac{\Cal{U}+J^{m-l}\Cal{V}}
{\Cal{Z}+J^{m-l}\Cal{W}},
$$
\noindent
where $\Cal{Z} \subseteq \Cal{U}$ and $\Cal{W} \subseteq \Cal{V}$ are finite, 
graded $S_1$-modules.

\endproclaim

\demo{Proof} It essentially goes as in Proposition 3 in \cite {ET}. \qed

\enddemo

\proclaim{Proposition 2} Let $R$, $S_1$, $J$ be as in Proposition 1.
Let $\Cal{T}$ be a graded $S_1$-module, and $\Cal{U}$, 
$\Cal{V}$, $\Cal{W}$, $\Cal{Z}$ be finite graded $S_1$-submodules of $\Cal{T}$.
Assume that $\Cal{Z} \subseteq \Cal{U}$, and that $\Cal{W} \subseteq \Cal{V}$, 
and denote

$$
{\Cal L}_m:= \frac{\Cal{U}+J^m\Cal{V}}{\Cal{Z}+J^m\Cal{W}}.
$$
Then, if ${({\Cal L}_m)}_n$, the $n$-th degree homogeneous component 
of ${\Cal L}_m$, has finite length for all large values of $m$ 
and $n$, $\l({({\Cal L}_m)}_n)$ is eventually given by a polynomial in 
$m$ and $n$.
\endproclaim

\demo{Proof} It follows the same path as Lemma 2, (b) in \cite{ET}. \qed
\enddemo

\proclaim{Proposition 3} Let $R$ be a Noetherian ring, $I$, $J \subseteq R$ 
ideals, $M$, $N$ be finite $R$-modules, and $i \geq 0$. If 
$\t_i(I^nM, N/J^mN)$ has finite length for all $m$, $n \gg 0$, then this 
length is eventually given by a polynomial in $m$, $n$.
\endproclaim

\demo{Proof} Take an $S_1$-graded resolution by finite free $S_1$-modules of 
the finite graded $S_1$-module $\oplus_{n=0}^{\infty} I^nM$. Tensor it with 
$N/J^mN$, in two steps, first with $N$, (call the resulting $S_1$-complex 
$\Cal C$), then with $R/J^m$. The part giving 
$\t_i^R(\oplus_{m=0}^{\infty} I^nM, N/J^mN)$, looks just like the situation 
described in Proposition 1. Therefore, by Proposition 1, we see that
$$
\t_i^{R}(\oplus_{n=0}^{\infty} I^nM, N/J^mN)=\frac{\Cal U + J^{m-l} \Cal V}
{\Cal Z + J^{m-l}\Cal W}, 
$$
\noindent
for some $l$, all $m \geq l$, where $\Cal U$, $\Cal V$, $\Cal Z$ and $\Cal W$ 
are all finite graded $S_1$-modules.
It follows that 

$$
\t_i^R(I^nM, N/J^mN)=\frac{{\Cal U}_n + J^{m-l} {\Cal V}_n}
{{\Cal Z}_n + J^{m-l}{\Cal W}_n},
$$
\noindent
by looking at homogeneous pieces of degree $n$ in the previous $\t$ formula. 
Thus, the conclusion follows from Proposition 2. \qed
\enddemo

\proclaim{Lemma 4} Let $(R,P)$ be Noetherian, local, $I$, $J \subset R$ 
ideals, $i \geq 0$. Then, for two finite $R$-modules $M$, $N$, we have:
\roster
\item "(a)" The image of the induced map 
$$
\t_i(I^nM,N) @>{H(f_i)}>> \t_i(M,N)
$$
\noindent
is of the form $I^{n-k}A$, for some $k \geq 0$ and $n \geq k$, where $A$ is 
the image of the map $\t_i(I^kM,N) @>{H(f_i)}>> \t_i(M,N)$ .

\item "(b)" The image of the induced map 
$$
\t_i(M,N) @>{H(g_i)}>> \t_i(M,N/J^mN)
$$
\noindent

has the form
$$
\frac{\t_i(M,N)+J^mB}{J^mB},
$$
\noindent

for some module $B$, such that $\t_i(M,N) \subseteq B$.
\endroster
\endproclaim
\demo{Proof} (a) Let 

$$
 \cdots @>>> R^{\beta_{i+1}} @>>> R^{\beta_i} @>>> R^{\beta_{i-1}} @>>> \cdots
\tag 1
$$
\noindent

be a free resolution of $N$. Then we have the following commutative diagram

$$
\CD
\cdots @>>> I^nM^{\beta_{i+1}} @>{\psi_n}>> I^nM^{\beta_i} @>{\phi_n}>> 
I^nM^{\beta_{i-1}} @>>> \cdots \\
  @.         @ VVV                   @ VV{f_i}V               @VVV 
  @.   \\
\cdots @>>> M^{\beta_{i+1}} @>{\psi}>> M^{\beta_i} @>{\phi}>> M^{\beta_{i-1}} @>>> \cdots .\\
\endCD
$$
\noindent
Let $K = \ker \phi$ and $L= \i \psi$, so $\t_i(M,N) = K/L$.
We also have that $\ker \phi_n = K \cap I^nM^{\beta_i}$ and $\i \psi_n = I^nL$, 
and thus $\t_i(I^nM,N)= (K \cap I^nM^{\beta_i})/I^nL$.
It follows that

$$
 \i (H(f_i))= \frac{K \cap I^nM^{\beta_i}+L}{L} = 
\frac{I^{n-k}(K \cap I^kM^{\beta_i})+L}{L}, 
$$
\noindent
for some $k$ and all $n \geq k$. 
Note that this is of the form $I^{n-k}A$, where $A$ is the image of the map 
$\t_i(I^kM,N) @>{H(f_i)}>> \t_i(M,N)$, as stated.

(b) Now assume that (1) gives a free resolution of $M$, and tensor it 
with $N/J^mN$. We get 

$$
\CD
\cdots @>>> N^{\beta_{i+1}} @>\psi>> N^{\beta_i} @>\phi>> N^{\beta_{i-1}} @>>> 
\cdots \\
@.             @VVV                   @VV{g_i}V           @VVV         @. \\
\cdots @>>> N^{\beta_{i+1}}/J^mN^{\beta_{i+1}} @>\psi_m>> 
N^{\beta_{i}}/J^mN^{\beta_i} @>\phi_m>> N^{\beta_{i-1}}/J^mN^{\beta_{i-1}} @>>> 
\cdots \\
\endCD
$$
\noindent
Again, if we denote $K=\ker \phi$ and $L= \i \psi$, then $\t_i(M,N) =K/L$ and, 
moreover, we obtain that 

$$
  \ker \phi_m = 
\frac{K+ J^{m-l}(\phi^{-1}(J^lN^{\beta_{i-1}}))}{J^mN^{\beta_i}},
$$
\noindent
for some $l$ and $m \geq l$.

We also get 
$$
  \i \psi_m= \frac{L+J^mN^{\beta_i}}{J^mN^{\beta_i}},
$$
\noindent
so 
$$
\t_i(M,N/J^mN)=\frac{K+ J^{m-l}(\phi^{-1}(J^lN^{\beta_{i-1}}))}
{L+J^mN^{\beta_i}}.
$$
\noindent
It follows that 

$$
\i H(g_i)=\frac{K+J^mN^{\beta_i}}{L+J^mN^{\beta_i}} \cong \frac{\t_i(M,N)+J^mB}
{J^mB},
$$
\noindent
where $B=N^{\beta_i}/L$. Of course, $\t_i(M,N) \subseteq B$. \qed
\enddemo

 The next Proposition is an extended version of the following well-known 
result:
Let $(R, P)$ be Noetherian, local, and $I \subseteq R$ an ideal. If $L$, $M$ 
are finitely generated modules, $L$ of finite length, then, for any 
$i \geq 0$, the natural map 
$\t_i(I^nM, L) @>>> \t_i(M, L)$ is zero, for $n \gg 0$ (see \cite{GL}).

\proclaim{Proposition 5} Let $(R,P)$ be a Noetherian, local ring. Let 
$I \subset R$ be an ideal, $M$, $N$ two finite $R$-modules and $i\geq 0$, 
fixed.
Then the following are equivalent:
\roster
\item "(a)"  $I \subseteq rad (ann_R \t_i(M,N))$.
\item "(b)"  $I \subseteq rad (ann_R \t_i(I^kM,N))$, for some $k \geq 0$.
\item "(c)"  $I \subseteq rad (ann_R \t_i(I^nM,N))$, for all $n \geq 0$.
\item "(d)"  $I \subseteq rad (ann_R \i (\t_i(I^nM,N) @>>> \t_i(M,N)))$, 
for all $n \geq 0$.
\item "(e)"  $\i (\t_i(I^nM,N) @>>> \t_i(M,N))=0$, for all $n \gg 0$.
\endroster
\endproclaim

\demo{Proof} Clearly, (c) implies (a) and (b). Conversely, consider the long 
exact sequence 
$$
\cdots @>>> \t_{i+1}(M/I^nM,N)@>\partial >> \t_i(I^nM,N) @>\alpha >> \t_i(M,N) 
@>\beta >> \t_i(M/I^nM,N) @>>> \cdots
$$
\noindent
(a) implies  (b), (c) follows by considering $\alpha$ and $\partial$, 
since $I \subseteq rad(ann_R \t_j(M/I^nM, N))$ for all $n \geq 0$. 
(b) implies (a) follows from (c) implies (a).
\newline
(a) implies (d) and (d) implies (a) are immediate, considering 
$\alpha$. 
\newline
(e) implies (a): if $\alpha=0$, then $\beta$ is an injection, so the 
conclusion follows.
\newline
(a) implies (e) follows from Lemma 4(a). \qed 
\enddemo

Here is the main result of this paper: 

\proclaim{Theorem 6} Let $(R, P)$ be Noetherian, local, $I$, $J \subseteq R$ 
two ideals, $M$, $N$ finitely generated $R$-modules, $i \geq 0$.  
 Assume that $ann_R M + ann_R N + J$ is $P$-primary. Then,
$$
   \l(\t_i(M/I^nM, N/J^mN))
$$
\noindent
is eventually given by a polynomial in $m$ and $n$ if and only if 
$I \subseteq rad(ann_R \t_j(M, N))$, for $j \in \{ i-1, i \}$.
\endproclaim

\demo{Proof} Consider the long exact sequence 
$$
  \cdots @ >>> \t_i(I^nM, N/J^mN) @ >\alpha_i^{m,n} >> \t_i(M, N/J^mN) 
@ >>> \t_i(M/I^nM, N/J^mN) @ >>> 
$$
$$
\t_{i-1}(I^nM ,N/J^mN) @ >\alpha_{i-1}^{m,n}>> \t_{i-1}(M, N/J^mN) @ >>> \cdots
$$

We already know that the lengths of the modules above, save the one in the 
middle, are (eventually) given by polynomials in one or two variables (see 
Proposition 3). Thus, we have 
$$
\align
  &\l(\t_i(M/I^nM, N/J^mN)) = \lbrack \l(\t_i(M, N/J^mN)) - 
\l(\i \alpha_i^{m,n}) \rbrack + \l(\ker \alpha_{i-1}^{m,n})\tag 2\\
&= \lbrack \l(\t_i(M, N/J^mN)) - \l(\i \alpha_i^{m,n})  \rbrack + 
\lbrack \l(\t_{i-1}(I^nM, N/J^mN)) - \l(\i \alpha_{i-1}^{m,n})  \rbrack.
\endalign
$$
Therefore, we need to examine $\l(\i \alpha_j^{m,n})$, for 
$j \in \{i-1, i \}$.
Consider the following commutative diagram:
$$
\CD
 \t_i(I^nM, N) @>\psi^{m,n}>> \t_i(I^nM, N/J^mN) @>\phi^{m,n}>> 
\t_{i-1}(I^nM, J^mN) \\
 @ VV{\sigma^{m,n}}V @ VV{\alpha_i^{m,n}}V @  VV{\tau^{m,n}}V \\
\t_i(M, N) @>\theta^{m,n}>> \t_i(M, N/J^mN) @>>> \t_{i-1}(M, J^mN) \\
 @.        @ VV{\pi_i^{m,n}}V @. \\
 @. \t_i(M, N/J^mN)/\alpha_i^{m,n}(L_{m,n}) @. \\
 @.        @ VVV @. \\
 @.        0                  @.
\endCD
\tag 3
$$
\noindent
where, $L_{m,n}=\i \psi^{m,n} = \ker \phi^{m,n}$. \newline

Note that the commutative diagram $(3)$ is a homogeneous piece of the diagram 
$(3^{\prime})$ below. That's because $\t_i^R$ is additive, and the natural maps 
in $(3)$ commute with the action of $I$ and $J$ on the modules occurring in 
this diagram. It follows that the diagram $(3^{\prime})$ is a commutative 
diagram of bigraded $S_2$-modules and maps.

$$
\CD
 \t_i^R(\Cal{I} \Cal{M}, \Cal{N}) 
@>\psi>> \t_i^R(\Cal{I}\Cal{M}, \Cal{N}/\Cal{J}\Cal{N}) 
@>\phi>> \t_{i-1}^R(\Cal{I}\Cal{M}, \Cal{J}\Cal{N}) \\
 @ VV{\sigma}V @ VV{\alpha_i}V @  VV{\tau}V \\
\t_i^R(\Cal{M}, \Cal{N}) @>\theta>> 
\t_i^R(\Cal{M}, \Cal{N}/\Cal{J}\Cal{N}) 
@>>> \t_{i-1}^R(\Cal{M}, \Cal{J}\Cal{N}) \\
@.        @ VV{\pi_i}V @. \\
 @. \t_i^R(\Cal{M}, \Cal{N}/\Cal{J}\Cal{N})/\alpha_i(\Cal{L}) @. \\
 @.        @ VVV @. \\
 @.        0                  @. ,
\endCD
\tag {$3^{\prime}$}
$$
where $\Cal{L}=\oplus_{m,n=0}^\infty L_{m,n}$.

Observe now that $\pi_i \circ \alpha_i$ factors through the image of $\phi$, 
which is a finitely generated, bigraded $S_2$-module (since 
$\t_{i-1}^R(\Cal{I}\Cal{M}, \Cal{J}\Cal{N})$ is so), hence 
$\i (\pi_i \circ \alpha_i)$ is a finite, bigraded $S_2$-module. Then 
$\l(\i (\pi_i \circ \alpha_i)^{m,n})$ is eventually given by a polynomial, by 
classical theory.

Note that $\l(\i \alpha_i^{m,n}) 
= \l(\i (\pi_i \circ \alpha_i)^{m,n}) + \l(\alpha_i^{m,n}(L_{m,n}))$, 
and a similar equality holds for $i-1$ in place of $i$. From (2) and what we 
have just seen, it follows that $\l(\t_i(M/I^nM, N/J^mN))$ is eventually 
given by a polynomial, if and only if the same is true of 
$\l(\alpha_{i-1}^{m,n}(L_{m,n}))  + \l(\alpha_i^{m,n}(L_{m,n}))$. 

We now examine $\l(\alpha_i^{m,n}(L_{m,n}))$. From $(3)$, we find that 
$$
\alpha_i^{m,n}(L_{m,n}) = \alpha_i^{m,n}(\psi^{m,n}(\t_i(I^nM, N))) = 
(\theta \circ \sigma)^{m,n}(\t_i(I^nM, N)).
\tag 4
$$
\noindent
From Lemma 4, (a) and (b), we get that 
$$
(\theta \circ \sigma)^{m,n}(\t_i(I^nM, N))=\frac{I^{n-k}A + J^mB}{J^mB}=
\frac{I^{n-k}A}{I^{n-k}A \cap J^mB},
\tag 5
$$

\noindent
for some $k \geq 0$ and $n \geq k$, where 
$A=\i (\t_i(I^kM, N) @ >>> \t_i(M, N))$.\newline 

We now claim that $\l(I^{n-k}A/I^{n-k}A \cap J^mB)$ is identically zero, 
for $m$, $n \gg 0$ if and only if it is polynomial for $m, n \gg 0$, 
if and only if $I \subseteq rad(ann_R \t_i(M, N))$. To prove this claim, 
assume $I \subseteq rad(ann_R \t_i(M, N))$. Then $I^{n-k}A=0$, for large $n$, 
and so $\l(I^{n-k}A/I^{n-k}A \cap J^mB)=0$, hence polynomial, for $n \gg 0$ 
and all $m$. It remains to check that, if 
$I \nsubseteq rad(ann_R \t_i(M, N))$, then $\l(I^{n-k}A/I^{n-k}A \cap J^mB)$ 
is nonzero and not given by a polynomial, for all $m, n \gg 0$. Indeed, by 
Proposition 5, (1) $\Leftrightarrow$ (3), we know that 
$I \nsubseteq rad (ann_R \i (\t_i(I^nM,N) @>>> \t_i(M,N)))$, 
for all $n$, so $I^{n-k}A \neq 0$ for all $n \geq k$. 

Now, since $ann_R M + ann_R N + J$ is $P$-primary, there is a $l \geq 0$, 
such that $I^l \subseteq ann_R M + ann_R N + J$. It follows that, for 
$n \geq lm + k$, we have 
$$
I^{n-k} \subseteq J^m + ann_R M + ann_R N,
$$
\noindent
so
$$
I^{n-k}A \subseteq J^mA \subseteq J^mB,
$$
\noindent
since we know that $A \subseteq B$.

Thus, for $n \geq lm + k$, $l$ and $k$ fixed, 
$\l(I^{n-k}A/I^{n-k}A \cap J^mB)$ vanishes. On the other hand, note that, 
for every $n \geq k$, $I^{n-k}A / I^{n-k}A \cap J^mB \neq 0$, for all 
$m \gg 0$. This is so since, for every 
$n \geq k$, $n$ fixed, $I^{n-k}A \cap J^mB \subsetneq I^{n-k}A$ for all 
large $m$, by Krull's Intersection Theorem. Hence 
$\l(I^{n-k}A/I^{n-k}A \cap J^mB) \neq 0$, for every $n \geq k$ and $m \gg 0$. 
This proves the claim, since we proved that, above the line $d: n = lm +k$ in 
the $(m, n)$-plane, $\l(I^{n-k}A/I^{n-k}A \cap J^mB)$ always vanishes, for 
large $m$ and $n$, while below this line, the length in question is nonzero, 
in case $I \nsubseteq rad(ann_R \t_i(M, N))$.

Finally, note that both terms of the form $\l(I^{n-k}A/I^{n-k}A \cap J^mB)$ 
occurring in the formula (2) of $\l(\t_i(M/I^nM, N/J^mN))$ (also see $(4)$ and 
$(5)$), actually occur with the same sign. By the claim, it follows 
that the sum of these two terms vanishes for all large $m$ and $n$, if 
$I \subseteq rad(ann_R \t_i(M, N))\cap rad(ann_R \t_{i-1}(M, N)) $. On the 
other hand, if 
$I \nsubseteq rad(ann_R \t_i(M, N))\cap rad(ann_R \t_{i-1}(M, N))$, then the 
sum in question vanishes above both lines $d: n = lm +k$, 
$d^\prime :n = l^\prime m +k^\prime $, (one line for each term), but it is 
nonzero below both these lines, $d$ and $d^\prime$. This means that 
$\l(\t_i(M/I^nM, N/J^mN))$ can only then be (eventually) 
polynomial, when both terms of the form $\l(I^{n-k}A/I^{n-k}A \cap J^mB)$ 
vanish. And this happens if and only if $I \subseteq rad(ann_R \t_j(M, N))$, 
for $j \in \{ i-1, i\}$, as stated. \qed   
\enddemo

The proof of Theorem 6 yields the following interesting corollary:

\proclaim{Corollary 7} Let $(R, P)$ be Noetherian, local, $I$ an ideal  
$M$, $N$ two finite $R$-modules and $i \geq 0$. Assume that 
$I + ann_R M + ann_R N$ is P-primary. Then 
$$
   \l(\t_i(M/I^nM, N/I^nN))
$$
\noindent
is given by a polynomial, for $n \gg 0$.
\endproclaim

\demo{Proof} Note that, by the proof of Theorem 6, we only have to look at 
each of the two (similar) terms in $\l(\t_i(M/I^nM, N/J^mN))$, that turned 
out  not to be polynomial, in general. If in each of them we set $J = I$ and 
$m = n$, we get two terms, each of which looks like 
 
$$
   \l(\frac {I^{n-k}A}{I^{n-k}A \cap I^nB}).
$$
\noindent
It is immediate, by the Artin-Rees Lemma, that 
$\oplus_{n=0}^{\infty} I^{n-k}A/I^{n-k}A \cap I^nB$ is a finite graded module 
over the Rees ring ${\Cal R}_I = \oplus_{n=0}^\infty I^n$, hence the 
conclusion. \qed 
\enddemo

\proclaim{Corollary 8} Assume that both $I + ann_R M + ann_R N$ and 
$J + ann_R M + ann_R N$ are $P$-primary, in the statement of Theorem 6. 
Then $\l(\t_i(M/I^nM, N/J^mN))$ is eventually given by a polynomial if and 
only if $\t_j(M, N)$ has finite length for both $j=i$, $j=i-1$.
\endproclaim

\demo{Proof} $\l(\t_i(M/I^nM, N/J^mN))$ is eventually given by a polynomial 
if and only if $I \subseteq rad (ann_R \t_j(M, N))$, for 
$j \in \{ i-1, i \}$,  if and only if $I + ann_R M + ann_R N \subseteq $ 
\newline 
$rad (ann_R \t_j(M, N))$, for $j \in \{ i-1, i \}$,   if and only if 
$\t_j(M, N)$, has finite length for both $j=i-1$ and $j=i$. \qed 
\enddemo

\subhead{Remark}
\endsubhead
 From this corollary alone we could construct numerous examples in which 
\linebreak
$\l(\t_i(M/I^nM, N/J^mN))$ is not eventually polynomial. It suffices to take 
$I$ and $J$ to be $P$-primary ideals and $M$, $N$ two finite $R$-modules with 
at least one of the two modules $\t_i(M,N)$ and $\t_{i-1}(M,N)$ not having 
finite length. Let us give two such examples of $\t_i(M/I^nM, N/J^mN)$
that have non-polynomial length, the second of which works for any
value of $i$.
 
 First, assume that $R$ has positive depth and dimension at least two. Take 
$x_1, x_2, \ldots , x_t$, $t \geq 1$ to be a regular sequence, such that the 
ideal generated by these elements is {\it not} $P$-primary. Take 
$M=R/(x_1, \ldots , x_t)^s$ and $N=R/(x_1, \ldots , x_t)^r$ for some 
$s \geq r \geq 1$. Then 
$\t_1(M,N)=(x_1, \ldots , x_t)^s/(x_1, \ldots , x_t)^{s+r}$ has finite
length if and only if $R/(x_1, \ldots , x_t)$ has finite length. This
is so because, by Rees' theorem, 
$(x_1, \ldots , x_t)^j/(x_1, \ldots , x_t)^{j+1}$ is a free 
$R/(x_1, \ldots , x_t)$-module, for all $j \geq 0$. Therefore 
$\t_1(M,N)$ can't have finite length by the choice of the regular sequence. 
Now take $I$ and $J$ any two $P$-primary ideals: by Corollary 8, the length 
of  $\t_i(M/I^nM, N/J^mN)$ is {\it not} given by a polynomial, for 
$i \in \{1, 2\}$.

 Secondly, assume that $R$ is neither regular, nor an isolated singularity. 
Then $R_Q$ is {\it not} regular for some non-maximal prime $Q$. Take
$M$ and $N$ to be any two finite $R$-modules, such that their annihilator is 
$Q$. Note that both $M_Q$ and $N_Q$ are direct sums of copies of the residue 
field of $R_Q$. Then $\t_i(M,N)$ cannot have finite length for any $i$. (For 
$i \geq 1$: this would imply that the localization at $Q$ of $\t_i(M,N)$ 
vanishes, giving that $R_Q$ is regular, contrary to the choice of $R$.) Now, 
Corollary 8 says that for any choice of two primary ideals $I$ and $J$, 
the length of $\t_i(M/I^nM, N/J^mN)$ is not polynomial for {\it all} 
$i\geq 0$. 

\proclaim{Theorem 9} Let $(R, P)$ be Noetherian local, $I$, $J \subseteq R$ 
ideals, $M$, $N$ finite $R$-modules and $i \geq 0$. Assume that $M \otimes N$ 
has finite length. Then 

$$
 \l(\t_i(M/I^nM, N/J^mN))
$$
\noindent
is given by a polynomial, for $m$, $n \gg0$.
\newline
Moreover, $\l(\t_i(M/I^nM, N/J^mN)) =$
$$
\l(\t_i(M, N))+\l(\t_{i-1}(I^nM, N))+
\l(\t_{i-1}(M, J^mN))+\l(\t_{i-2}(I^nM, J^mN)) 
$$
\endproclaim

\demo{Proof} The first statement follows immediately from Theorem 6, since, 
trivially, its hypotheses are met. For the last statement, let's observe 
that, there is a $k \geq 0$, such that, for all $m \geq 0$, and $n\geq k$,  
$\sigma^{m,n}$ in (3) is the zero map, by Proposition 5. 
It follows that $\alpha_i^{m,n}(\i \psi^{m,n})
=\alpha_i^{m,n}(\ker \phi^{m,n})=0$, hence $\alpha_i^{m,n}$ factors through 
$\i \phi^{m,n}$, and thus (as before) $\l(\i \alpha_i^{m,n})$ is eventually 
given by a polynomial in $m$, $n$. Finally, by Proposition 5 again, we see 
that for each {\it fixed} $m$, $\i (\alpha_i^{m,n})$ vanishes for $n \gg 0$. 
Therefore, $\i (\alpha_i^{m,n})$ is identically zero, for all large $m$ and 
$n$. 

We also have the long exact sequence 
$$
\cdots @>>> \t_i(I^nM, N/J^mN) @>{\alpha_i^{m,n}}>> \t_i(M, N/J^mN) @>>> 
\t_i(M/I^nM, N/J^mN) @>>> 
$$
$$
\t_{i-1}(I^nM, N/J^mN) @>{\alpha_{i-1}^{m,n}}>> \t_{i-1}(M, N/J^mN) @>>> 
\cdots ,
$$ 
\noindent
and we now know that $\alpha_i^{m,n}=\alpha_{i-1}^{m,n}=0$ for $m$, 
$n \gg 0$. Then, 
$$
 \l(\t_i(M/I^nM, N/J^mN))=\l(\t_i(M, N/J^mN))+\l(\t_{i-1}(I^nM, N/J^mN)).
\tag 6
$$
\noindent
We apply this trick two more times. We have 
$$
 \dots @>>> \t_i(M, J^mN) @>0>> \t_i(M, N) @>>> \t_i(M, N/J^mN) @>>> 
$$
$$
\t_{i-1}(M, J^mN) @>0>> \t_{i-1}(M, N) @>>> \cdots ,
\tag 7
$$
\noindent
where the maps marked as 0 are so by Proposition 5. We get that 
$$
\l(\t_i(M, N/J^mN))=\l(\t_i(M, N))+\l(\t_{i-1}(M, J^mN)).
\tag 8
$$
\noindent
Replacing $M$ by $I^nM$ in $(7)$ and using the fact that 
$\oplus_{m,n=0}^{\infty} \t_i(I^nM, J^mN)$ is a finite bigraded $S_2$-module, 
we see that the maps marked as 0 will remain so, for every $n$ and large $m$, 
again by Proposition 5. 
We then get that 
$$
\l(\t_{i-1}(I^nM, N/J^mN))=\l(\t_{i-1}(I^nM, N))+\l(\t_{i-2}(I^nM, J^mN)).
\tag 9
$$
\noindent
Putting together (6), (8) and (9), we obtain 
$$
\l(\t_i(M/I^nM, N/J^mN)) 
$$
$$
= \l(\t_i(M, N))+\l(\t_{i-1}(I^nM, N))+
\l(\t_{i-1}(M, J^mN))+\l(\t_{i-2}(I^nM, J^mN)),
$$
\noindent
as stated. \qed 
\enddemo

Note that this also yields a direct proof of the first statement of this 
theorem, since the four terms on the right-hand side of the equality above 
are eventually given by polynomials, by classical theory of finite (bi)graded 
modules.

Finally, we give an upper bound for the degree of the polynomial that arises 
in Corollary 8. Note that this estimate also applies to the case of Theorem 9.

\proclaim{Proposition 10} Assume the hypotheses in Corollary 8 and
suppose that the length of $\t_i(M/I^nM, J^mN)$ is given by a polynomial, for 
$m$, $n \gg 0$. Then 

$$
   \deg \l(\t_i(M/I^nM, J^mN)) \leq \ell_M(I) + \ell_N(J) -2.
$$
\endproclaim

\demo{Proof} This is a rather crude estimate, based on the one-variable case. 
We simply apply Corollary 4 in \cite {ET}, separately, for fixed, large 
enough values of $m$ and $n$, then add. For the exact degree in some special 
cases (in one variable, though), see \cite{TM}. \qed 
\enddemo

\head{Acknowledgements} 
\endhead

 The work in this paper forms part of 
my dissertation written at the University of Kansas. 
I would particularly like to thank my advisor, Prof. D. Katz, for his 
constant patience and steady encouragements. I would also like to thank both 
D. Katz and A. Vraciu for their comments on this paper.

 Special thanks also go to the referee for a number of valuable suggestions.

\enddocument